\documentclass{elsart}               % The use of LaTeX2e is preferred.

\usepackage{graphicx}          % Include this line if your
\usepackage{amsmath}
\usepackage{amssymb}

\begin{document}

\begin{frontmatter}

\title{Practical and efficient method for fractional-order unstable pole-zero cancellation in linear feedback systems \thanksref{footnoteinfo}} % Title, preferably not more

\thanks[footnoteinfo]{Corresponding author F. Merrikh-Bayat. Tel.
+98-241-5154061. Fax +98-241-5152762.}

\author[Paestum]{Farshad Merrikh-Bayat}\ead{f.bayat@znu.ac.ir}    % Add the
%\author[Rome]{Julius Caesar}\ead{julius@caesar.ir},               % e-mail address
%\author[Baiae]{Publius Maro Vergilius}\ead{vergilius@culture.ir}  % (ead) as shown

\address[Paestum]{Department of Electrical and Computer Engineering, University of Zanjan, Zanjan, IRAN}  % Please supply
%\address[Rome]{Senate House, Rome}             % full addresses
%\address[Baiae]{The White House, Baiae}        % here.

                             % keyword list or with the
                                          % help of the Automatica
                                          % keyword wizard

\begin{abstract}                          % Abstract of not more than 200 words.
As a very well-known classical fact, non-minimum phase zeros of
the process put some limitations on the performance of the
feedback system. The source of these limitations is that such
non-minimum phase zeros cannot be cancelled by unstable poles of
the controller since such a cancellation leads to internal
instability. The aim of this paper is to propose a method for
fractional-order cancellation of non-minimum phase zeros of the
process and studying its properties. It is specially shown that
the proposed cancellation strategy increases the phase and the
gain margin without leading to internal instability. Since the
systems with higher gain and phase margin are easier to control,
the proposed method can be used to arrive at more effective
controls.
\end{abstract}

\end{frontmatter}

\section{Introduction}
Non-minimum phase zeros appear unavoidably in some important
processes such as steam generators \cite{astrom}, aircrafts
\cite{hauser}, flexible-link manipulators \cite{kwon}, and
continuous stirred tank reactors \cite{kravaris}. As a very
well-known classical fact, non-minimum phase zeros of the process
limit the performance of the feedback system in different ways
\cite{seron1}-\cite{middleton}. For instance, these limitations
can be concluded from the classical root-locus method
\cite{hoagg}, asymptotic LQG theory \cite{qiu}, and waterbed
effect phenomena \cite{doyle}. In the field of linear
time-invariant (LTI) systems, the source of these limitations is
that the non-minimum phase zero of the process cannot be cancelled
by unstable pole of the controller since such a cancellation leads
to internal instability \cite{kailath}.

%\cite{seron1},\cite{seron2},\cite{qiu}, \cite{middleton}
% \cite{fliness}
So far, various methods have been developed for the control of
processes with non-minimum phase zeros (see, for example,
\cite{Aguiar}-\cite{hiddabi} and the references therein for more
information). Clearly, according to the high achievement of
feedback control systems, it is strictly preferred to develop more
effective methods to the control of non-minimum phase systems
based on the feedback strategy. The aim of this paper is to
propose a modified feedback control strategy for non-minimum phase
processes, which is based on subjecting the non-minimum phase zero
of the process to a kind of cancellation. More precisely, it will
be shown that the non-minimum phase zero (unstable pole) of the
process can partly be cancelled by fractional-order pole (zero) of
the controller without leading to internal instability. It will
also be shown that the non-minimum phase zero of the process can
be cancelled to an arbitrary degree by the pole of the
fractional-order controller, only at the cost of using a more
complicated setup. Interesting observation, which is also
supported by mathematical discussions, is the fractional-order
cancellation of the non-minimum phase zero can considerably
increase the phase and gain margin, and consequently, make the
system easier to control.

The rest of this paper is organized as follows. Section
\ref{sec_main} contains the main results of paper. Proposed method
for fractional-order cancellation of the non-minimum phase zero is
presented in this section and it is shown that the proposed
cancellation strategy can improve the robustness of the feedback
system. Effect of this method on time-domain undershoots in the
step response of both open-loop and closed-loop systems is also
discussed in Section \ref{sec_main}. Three illustrative examples
are studied in Section \ref{sec_exam}, and finally, Section
\ref{sec_conc} concludes the paper.

\section{Main results}\label{sec_main}
Consider a LTI process with transfer function $P(s)$, input $u(t)$
and output $y(t)$. Suppose that $P(s)$ has a positive real zero of
order one at $s=\lambda$, that is $P(\lambda)=0$ and
$P'(\lambda)\neq 0$ where $\lambda$ is a positive real number.
Such a transfer function can be decomposed as
\begin{equation}\label{proc}
P(s) =\left(1-\frac{s}{\lambda}\right)\widetilde{P}(s).
\end{equation}
In the above equation the term $1-s/\lambda$ can be expanded using
fractional powers of $s$ in infinite many different ways. A
straightforward approach is to write it as
\begin{align*}
1-\frac{s}{\lambda}&=\left[1-\left(\frac{s}{\lambda}
\right)^{1/2}\right] \left[1+\left(\frac{s}{\lambda}
\right)^{1/2}\right]\\&=\left[1-\left(\frac{s}{\lambda}
\right)^{1/4}\right]\left[1+\left(\frac{s}{\lambda}
\right)^{1/4}\right] \left[1+\left(\frac{s}{\lambda}
\right)^{1/2}\right]\\&=\left[1-\left(\frac{s}{\lambda}
\right)^{1/8}\right]\left[1+\left(\frac{s}{\lambda}
\right)^{1/8}\right]\left[1+\left(\frac{s}{\lambda}
\right)^{1/4}\right] \left[1+\left(\frac{s}{\lambda}
\right)^{1/2}\right]\\&=\ldots,
\end{align*}
which yields
\begin{equation}\label{expand1}
1-\frac{s}{\lambda}=\left[1-\left(\frac{s}{\lambda}\right)
^{1/v}\right]\prod_{k=0}^{\log_2 (v/2)}\left [1+\left
(\frac{s}{\lambda}\right) ^{2^k/v}\right],
\end{equation}
where $\log_2 (v/2)$ is the base 2 logarithm of $v/2$, and $v$ can
be considered equal to any number in the form of $v=2^h$,
$h=1,2,3,\ldots$. The expression in the right-hand side of
(\ref{expand1}) has exactly $1+\log_2(v)$ roots distributed on a
Riemann surface with $v$ Riemann sheets, where the origin is a
branch point of order $v-1$ \cite{silverman}. Note that among
these roots only the root of $1-(s/\lambda)^{1/v}=0$ is located on
the first Rimann sheet and other roots are located on other sheets
\cite{silverman}.

%It is a well-known fact in the field of fractional-order systems
%that the features of system are mainly determined by the poles and
%zeros of

Substitution of (\ref{expand1}) in (\ref{proc}) and dividing both
sides of the resulted equation to
\begin{equation}
Q_{\lambda,v}(s) \triangleq
\prod_{k=0}^{\log_2(v/2)}\left[1+\left(\frac{s}{\lambda}\right)
^{2^k/v}\right],
\end{equation}
yields
\begin{equation}\label{pf}
P_{f}(s)\triangleq \frac{P(s)}{Q_{\lambda,v}(s)}=
\left[1-\left(\frac{s} {\lambda}\right)^{1/v}\right]
\widetilde{P}(s).
\end{equation}
Note that $P(s)$ and $P_f(s)$ are exactly the same (in the sense
that they have the same poles and zeros and DC gains) except that
$P_f(s)$ has a \emph{weaker} non-minimum phase zero at
$s=\lambda$. In fact, the zero of $P_f(s)$ at $s=\lambda$ is
weaker than the zero of $P(s)$ at $s=\lambda$ since it makes the
system less non-minimum phase (see the discussions below). In the
rest of this paper when the process transfer function is applied
in series with a system with transfer function
$1/Q_{\lambda,v}(s)$ we say that the process is subjected to a
fractional-order pole-zero cancellation.

In the following we discuss on the effects of fractional-order
cancellation of the non-minimum phase zero of the process on the
time and frequency domain characteristics of both open-loop and
closed-loop systems.

%First, they make the phase of the open-loop transfer function more
%negative by injecting a negative phase to the loop, and second,
%they increase the cutoff frequency of the open-loop transfer
%function by increasing its magnitude.

\subsection{Effect of the fractional-order unstable pole-zero cancellation on phase
margin}\label{subsec_pm} For two main reasons non-minimum phase
zeros of the process put a limitation on the robust stability of
the feedback system: First, they push the Bode phase plot of the
open-loop transfer function downward by injecting a negative phase
to the loop, and second, they increase the crossover frequency of
the Bode magnitude plot of the open-loop transfer function by
increasing its magnitude at all frequencies. These two reasons
together decrease the phase and the gain margin of the feedback
system. In the following, we show that applying the proposed
fractional-order cancellation strategy to the non-minimum phase
zero of the process can increase the phase margin of the feedback
system by partly removing both of the above-mentioned reasons.

Without any loss of generality consider the following processes:
\begin{equation}\label{p1}
P_1(s)=K\frac{1}{(1+s/\gamma_1)(1+s/\gamma_2)},
\end{equation}
\begin{equation}\label{p2}
P_2(s)=K\frac{1-s/\lambda}{(1+s/\gamma_1)(1+s/\gamma_2)},
\end{equation}
\begin{equation}\label{p3}
P_3(s)=K\frac{1-(s/\lambda)^{1/v}}{(1+s/\gamma_1)(1+s/\gamma_2)},
\end{equation}
where $K$, $\gamma_1$, $\gamma_2$, and $\lambda$ are positive real
constants such that $\gamma_1<< \lambda << \gamma_2$, and $v$ is
an integer constant greater than unity. As it can be observed,
$P_1(s)$ is minimum phase, $P_2(s)$ has a non-minimum phase zero
at $s=\lambda$, and $P_3(s)$ has a weaker non-minimum phase zero
(compared to the non-minimum phase zero of $P_2(s)$) at
$s=\lambda$ (note that $P_3(s)$ is obtained by applying the
proposed fractional-order unstable pole-zero cancellation to
$P_2(s)$).

\begin{figure}
\begin{center}
\includegraphics[width=8.4cm]{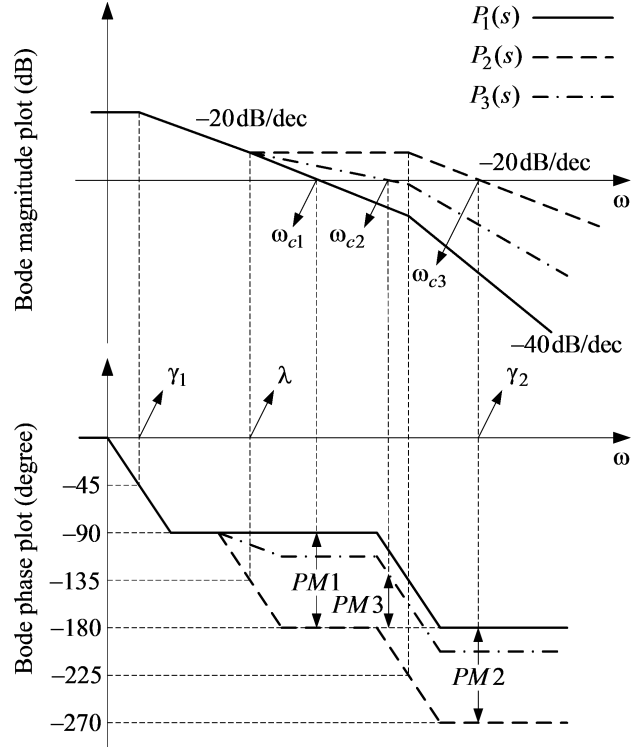}
\caption{The Bode phase and magnitude plot of the transfer
functions given in (\ref{p1})-(\ref{p3})}\label{fig_bode}
\end{center}
\end{figure}

Effect of the proposed cancellation strategy on the phase margin
can be deduced from the Bode plots, which are asymptotically
depicted in Figure \ref{fig_bode}. In this figure $PM1$, $PM2$,
and $PM3$ denote the phase margins corresponding to $P_1(s)$,
$P_2(s)$, and $P_3(s)$, respectively. Figure \ref{fig_bode}
clearly shows that at frequencies smaller than $\lambda$ all Bode
plots are almost the same but at frequencies larger than $\lambda$
the Bode magnitude plot of $P_3(s)$ decays faster than the Bode
magnitude plot of $P_2(s)$ while its Bode phase plot decays slower
than the Bode phase plot of $P_2(s)$. It turns out that
fractional-order cancellation of the non-minimum phase zero will
increase the phase margin and simultaneously decrease the gain
crossover frequency (note that in Fig. \ref{fig_bode} we have
$PM3>PM2$ and $\omega_{c2} <\omega_{c3}$). In the same manner it
can be easily verified that fractional-order cancellation of the
non-minimum phase zero will also increase the gain margin.

As we know, in dealing with feedback systems non-minimum phase
zeros of the process put a limitation on the gain of the
controller (recall from the root-locus method that closed-loop
poles move toward the open-loop zeros as the gain is increased).
On the other hand, the above discussions showed that the gain
margin of the feedback system is increased by applying the
proposed fractional-order cancellation to the non-minimum phase
zero of the process. Hence, it is possible to use controllers with
larger gains in the loop when such a cancellation is applied. It
turns out that the proposed fractional-order cancellation strategy
can partly remove one of the limitations put on the performance of
the feedback system by non-minimum phase zeros of the process.
Finally, note that although the results of this section are
obtained based on studying the Bode plots of two special transfer
functions, these are quite general observations and can be easily
extended to other cases.

\subsection{Effect of the fractional-order unstable pole-zero cancellation on the time-domain undershoot}\label{subsec_under}
In this section we show that applying the proposed
fractional-order unstable pole-zero cancellation to a non-minimum
phase process will reduce the time-domain undershoots in its step
response. Note that it can be easily concluded from the previous
discussions that the overshoots in the step response of the
closed-loop system will be decreased by subjecting the non-minimum
phase zero of the process to the proposed fractional-order
pole-zero cancellation. More precisely, in Section \ref{subsec_pm}
we showed that the proposed method for fractional-order
cancellation of non-minimum phase zeros has the property that
always increases the phase margin, and as we know increasing the
phase margin naturally decreases the time-domain overshoots.

Studying the effect of the proposed cancellation strategy on
time-domain undershoots is more complicated. Consider again the
process with transfer function $P(s)$ as given in (\ref{proc}).
Assuming that the system input, $u(t)$, is equal to the step
function and the system is initially at rest, the relative
undershoot $r_{us}$ is defined as \cite{lau}:
\begin{equation}\label{us_def}
r_{us}(y(.))\triangleq -\inf_{t\in (0,\infty)}
\left\{\frac{y(t)}{\overline{y}}\right\},
\end{equation}
where $\overline{y}$ is equal to the steady-state value of $y(t)$.
Clearly, the above definition provides us with a reasonable
criterion to compare the undershoot of different systems. The
Laplace transforms of $y(t)$ and $u(t)$, respectively denoted as
$Y(s)$ and $U(s)$, satisfy the following relation with the process
transfer function:
\begin{equation}\label{temp1}
Y(s)=\int_{t=0}^\infty y(t)e^{-st}dt=P(s)U(s).
\end{equation}
Substitution of $s=\lambda$ in (\ref{temp1}) and considering the
fact that $P(s)$ has a zero at $s=\lambda$ yields:
\begin{equation}\label{temp2}
\int_{t=0}^\infty y(t) e^{-\lambda t}dt=0.
\end{equation}
Suppose that $y(t)$ reaches its steady-state value at $t=T$, that
is $y(t)=\overline{y}$ for $t>T$ (clearly, the step response of a
LTI system cannot reach its steady-state value at a finite time,
but in practice $T$ can be approximated with the settling time of
the system step response). Using this assumption equation
(\ref{temp2}) can be written as
\begin{equation}
\int_{0}^Ty(t)e^{-\lambda
t}dt+\int_{T}^\infty\overline{y}e^{-\lambda t}dt=0,
\end{equation}
or equivalently,
\begin{equation}\label{temp3}
\int_0^T-\frac{y}{\overline{y}}~e^{-\lambda t}dt=\int_T^\infty
e^{-\lambda t}dt=\frac{e^{-\lambda T}}{\lambda}.
\end{equation}
Now, according to (\ref{us_def}) equation (\ref{temp3}) yields
\begin{equation}
\int_0^Tr_{us}e^{-\lambda t}dt\ge \frac{e^{-\lambda T}}{\lambda},
\end{equation}
which tunes out
\begin{equation}\label{temp4}
r_{us}\ge\frac{1}{e^{\lambda T}-1}.
\end{equation}
Equation (\ref{temp4}) concludes that the lower bound of the
relative undershoot is decreased by increasing the frequency of
the non-minimum phase zero, $\lambda$, and/or the settling time,
$T$. Both $P(s)$ and $P_f(s)$ (as defined in (\ref{pf})) have a
zero at $s=\lambda$ but it is expected that the settling time of
the (step) response of a system with transfer function $P_f(s)$ be
larger than the settling time of the (step) response of a system
with transfer function $P(s)$. This statement can easily be
concluded from the fact that the settling time of the step
response of a system is decreased by increasing its gain crossover
frequency (or equivalently, its bandwidth). Since the gain
crossover frequency of $P_f(s)$ is always smaller than the gain
crossover frequency of $P(s)$ it is expected that the response of
a system with transfer function $P(s)$ settles down faster than
the response of a system with transfer function $P_f(s)$. Hence,
according to (\ref{temp4}) the \emph{lower bound} on the relative
undershoot is decreased by subjecting the non-minimum phase zero
of the process to the proposed fractional-order cancellation. Note
that although decreasing the lower bound of undershoot does not
necessarily mean that the undershoot itself will also be
decreased, it shows the potential of the proposed approach for
this purpose.

The above result can be extended to feedback systems. The key idea
is to note that the bandwidth of the closed-loop system is
directly proportional to the gain crossover frequency of the
open-loop transfer function. Since applying the fractional-order
unstable pole-zero cancellation to the given non-minimum phase
process located in a feedback system decreases the gain crossover
frequency of the open-loop transfer function, we can conclude that
the bandwidth of the closed-loop system is also decreased by
performing such a cancellation. As a result, decreasing the
bandwidth of the closed-loop system leads to increasing the
settling time of the system response to step command and
consequently decreasing the lower bound on the relative
undershoot.

\subsection{Internal-stability analysis of the feedback system containing fractional-order unstable pole-zero cancellation}
Previous discussions showed that the proposed strategy for
fractional-order cancellation of the non-minimum phase zero of the
process can affect some important characteristics of the feedback
system. Specially, it was shown that such a cancellation increases
the phase and gain margin (and consequently, improves the
robustness of system) and potentially can decrease the undershoots
in time-domain response. In this section we show that the proposed
fractional-order unstable pole-zero cancellation method has the
advantage that can be used in feedback systems without leading to
internal instability.

\begin{figure}
\begin{center}
\includegraphics[width=8.4cm]{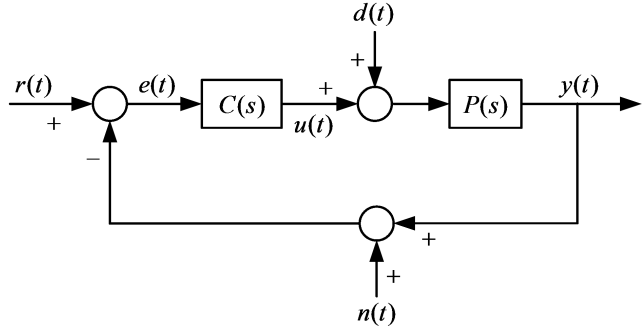}
\caption{The feedback system under
consideration}\label{fig_system}
\end{center}
\end{figure}

Consider the feedback system shown in Fig. \ref{fig_system} where
$P(s)$ is the process transfer function defined in (\ref{proc})
and $C(s)$ is the controller. Clearly, the system is internally
unstable if $C(s)$ has a pole at $s=\lambda$ (i.e., an unstable
pole-zero cancellation occurs in the loop). This fact can be
concluded by calculating, e.g., the transfer function from $r$ to
$u$, which definitely has an unstable pole at $s=\lambda$. Now,
assume that instead of the term $1-s/\lambda$ the denominator of
$C(s)$ contains the term $Q_{\lambda,v}(s) =
\prod_{k=0}^{\log_2(v/2)}\left[1+(s/\lambda)^{2^k/v}\right]$,
which is equivalent to a fractional-order pole-zero cancellation
between $C(s)$ and $P(s)$. In this case, instead of the pole at
$s=\lambda$ the transfer function from $r$ to $u$ has poles at the
roots of $Q_{\lambda,v}(s) =0$, which are distributed on a Riemann
surface. Important point that should be noted here is that the
stability of a system with fractional-order characteristic
equation $Q_{\lambda,v}(s) =0$ cannot be studied simply by
investigating the roots of this equation at the right half-plane
of the complex $s$-plane. In fact, according to the stability test
of Matignon \cite{matignon1} a system with fractional-order
characteristic equation
\begin{equation}\label{char1}
Q_{\lambda,v}(s) = \prod_{k=0}^ {\log_2(v/2)}
\left[1+\left(\frac{s}{\lambda}\right) ^{2^k/v}\right]=0,
\end{equation}
is stable if and only if all roots of the equation
$\widetilde{Q}_{\lambda,v}(w)=Q_{\lambda,v}(s)|_{s^{1/v}=w}=0$ lie
in the angular sector defined by
\begin{equation}\label{stab_sec}
|\arg (w)|>\frac{\pi}{2v},
\end{equation}
 in $w$-plane (in other words,
$Q_{\lambda,v}(s)=0$ must not have any poles in the closed right
half-plane of the first Riemann sheet). It can be easily verified
that the roots of $\widetilde{Q}_{\lambda,v}(w)=0$ are as follows:
\begin{equation}
w=e^{j\pi/2^k}\lambda^{1/v},\quad k=1,2,\ldots,\log_2(v/2).
\end{equation}
Since all of these roots are located in the sector of stability
defined in (\ref{stab_sec}) we can conclude that the proposed
fractional-order pole-zero cancellation method does not change the
status of internal stability of the feedback system (stability of
all other possible transfer functions in Fig. \ref{fig_system} can
be concluded in the same manner). For example, consider the
feedback system shown in Fig. \ref{fig_system} and suppose that
\begin{equation}
P(s)=\frac{s-1}{s+2}.
\end{equation}
Clearly, in this example $C(s)=1/(s-1)$ leads to internal
instability, while assuming
\begin{equation}
C(s)=\frac{1}{Q_{1,2}(s)}=\frac{1}{s^{1/2}+1},
\end{equation}
the transfer function from the reference input to control is
obtained as
\begin{equation}\label{temp5}
\frac{U(s)}{R(s)}=\frac{s+2}{(s^{1/2}+1)(s+s^{1/2}+1)},
\end{equation}
which is stable (note that substitution of $s^{1/2}=w$ in the
denominator of (\ref{temp5}) leads to the characteristic equation
$(w+1)(w^2+w+1)=0$, all roots of which are located in the sector
of stability defined as $|\arg(w)|>\pi/4$).

\begin{figure}
\begin{center}
\includegraphics[width=8.5cm]{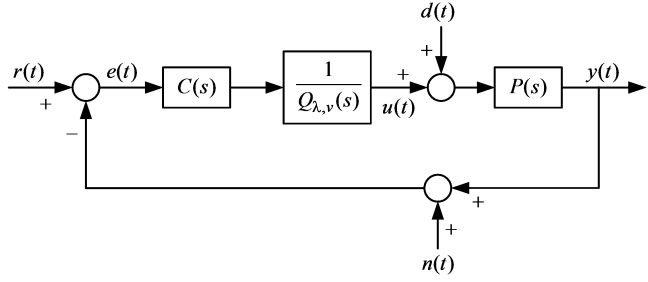}
\caption{Proposed structure for fractional-order cancellation of
the non-minimum phase zero of the process}\label{fig_prop}
\end{center}
\end{figure}

According to the previous discussions, the structure shown in Fig.
\ref{fig_prop} can be used for fractional-order unstable zero-pole
cancellation and feedback control of the non-minimum phase process
with transfer function $P(s)$, which has a non-minimum phase zero
at $s=\lambda$. In this structure, the final controller is equal
to the series combination of $C(s)$ and $1/Q_{\lambda,v}(s)$. In
the proposed method first  $Q_{\lambda,v}(s)$ should be determined
simply by assigning a suitable integer to $v$ $(v\ge 2)$. The
results obtained in Section \ref{subsec_pm} (and presented in Fig.
\ref{fig_bode}) show that the phase-lag of the transfer function
$P(s)/Q_{\lambda,v}(s)$ is decreased by increasing the value of
$v$, that is the transfer function $P(s)/Q_{\lambda,v}(s)$ becomes
less non-minimum phase (and consequently, an easier problem to
control) as the value assigned to $v$ is increased. However,
simulation results show that in many cases using $v=2$ or $v=4$
leads to satisfactory results.

After determining $Q_{\lambda,v}(s)$ the controller $C(s)$ can be
designed using any classical method assuming that the transfer
function of process is equal to $P(s)/Q_{\lambda,v}(s)$. Note that
since $P(s)/Q_{\lambda,v}(s)$ is a fractional-order transfer
function, the controller can be designed in two different ways:
One can apply an order-reduction algorithm to
$P(s)/Q_{\lambda,v}(s)$ to arrive at an  integer-order approximate
transfer function and then use a classical controller design
algorithm, or one can directly use the methods available to design
a controller for the given fractional-order process. Note also
that since $P(s)/Q_{\lambda,v}(s)$ has a smaller bandwidth
compared to $P(s)$, the controller designed for
$P(s)/Q_{\lambda,v}(s)$ naturally applies more control effort
compared to the similar controller designed for $P(s)$. A
reasonable approach to remove this difficulty is to apply the
proposed fractional-order pole-zero cancellation method to both
the non-minimum phase zero and a pole of process (an example of
this type is studied in Example 3 of Section \ref{sec_exam}).
Finally, note that when $P(s)$ has multiple non-minimum phase
zeros, say at $\lambda_1, \lambda_2,\ldots, \lambda_n$, it is
possible to consider a separate fractional-order pole-zero
canceller for each non-minimum phase zero.

%Since $P(s)/Q_{\lambda,v}(s)$ is less non-minimum phase compared
%to $P(s)$ it is expected

\subsection{Notes on realization of the proposed fractional-order pole-zero canceller}
One effective approach to realize the fractional-order pole-zero
canceller $Q_{\lambda,v}(s)$ is to approximate it with an
integer-order transfer function using an order reduction algorithm
(the Matlab function \textsf{invfreqs} is used for this purpose in
this paper). Note that although the approximate transfer function
obtained by using this method does not guarantee the (partial)
cancellation of any non-minimum phase zero or unstable pole of the
process, it can be used and will work in practice. The reason for
this statement is that the proposed fractional-order pole-zero
canceller actually affects the frequency response of the process
in a certain manner, and consequently, any other system that
applies the similar effects can be used instead. Hence, in
practice, we can simply substitute $Q_{\lambda,v}(s)$ with an
integer-order transfer function which has almost the same
frequency response as $Q_{\lambda,v}(s)$ in the bandwidth of the
process.

\section{Illustrative examples}\label{sec_exam}
Three illustrative examples are presented in this section in order
to confirm the results obtained in previous section.

%, as it was also mentioned in Section \ref{subsec_under},

\begin{figure}
\begin{center}
%[height=4cm]
\includegraphics[width=8.4cm]{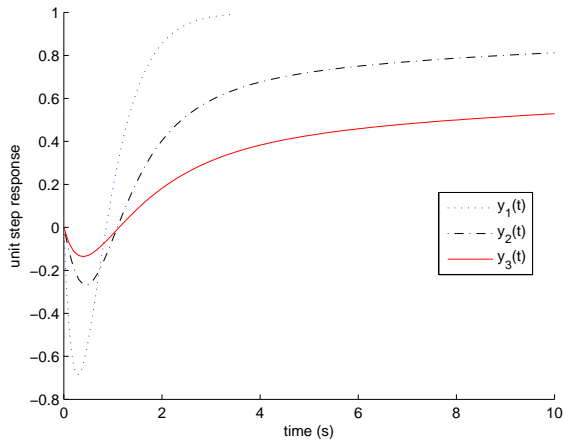}    % The printed column
\caption{Unit step response of three processes with transfer
functions given in (\ref{ex1_p1})-(\ref{ex1_p3}), corresponding to
Example 1}
\label{fig_ex1}                                 % Size the figures
\end{center}                                 % accordingly.
\end{figure}

\textbf{Example 1.} The aim of this example is to show that
fractional-order cancellation of the non-minimum phase zero of a
process leads to decreasing the undershoot and increasing the
settling time of its step response. For this purpose consider
three processes with the following transfer functions
\begin{equation}\label{ex1_p1}
P_1(s)=\frac{1-s}{(1+s/2)(1+s/3)},
\end{equation}
\begin{equation}
P_2(s)=\frac{P_1(s)}{1+s^{1/2}}=\frac{1-s^{1/2}}{(1+s/2)(1+s/3)},
\end{equation}
\begin{equation}\label{ex1_p3}
P_3(s)=\frac{P_1(s)}{(1+s^{1/2})(1+s^{1/4})}=\frac{1-s^{1/4}}{(1+s/2)(1+s/3)}.
\end{equation}
The corresponding unit step responses are shown by $y_1(t)$,
$y_2(t)$, and $y_3(t)$ in Fig. \ref{fig_ex1}. As it can also be
concluded from (\ref{temp4}), it is observed that applying the
proposed fractional-order pole-zero cancellation has decreased the
undershoot in the unit step response at the cost of increasing the
settling time. Since the zero of $P_3(s)$ is less non-minimum
phase compared to the zero of $P_2(s)$, $y_3(t)$ exhibits smaller
undershoot and larger rise time compared to $y_2(t)$.

%
%\begin{equation}
%P_2(s)=\frac{P_1(s)}{1-s}=\frac{1}{(1+s/2)(1+s/3)},
%\end{equation}

\begin{figure}
\begin{center}
%[height=4cm]
\includegraphics[width=8.4cm]{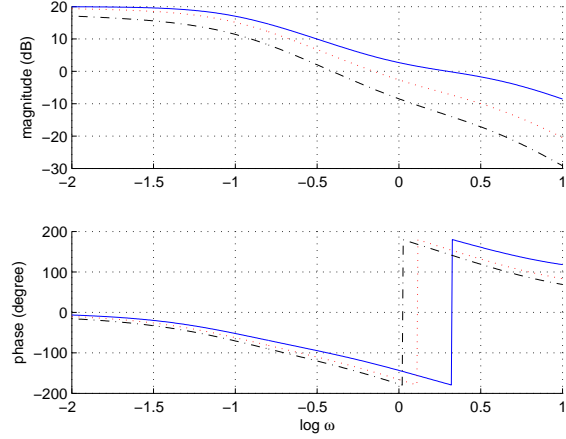}    % The printed column
\caption{Bode plot of $P(s)$ (solid curve), $P(s)/Q_{1,2}(s)$
(dotted curve), and $P(s)/Q_{1,4}(s)$ (dash-dotted curve),
corresponding to Example 2.}
\label{fig_ex2}                                 % Size the figures
\end{center}                                 % accordingly.
\end{figure}

\textbf{Example 2.} In this example we show that applying the
proposed fractional-order pole-zero cancellation strategy to a
non-minimum phase process can improve the robust stability of the
corresponding feedback system, and consequently make the control
problem easier. Consider the feedback system shown in Fig.
\ref{fig_prop} and assume that
\begin{equation}\label{p}
P(s)=\frac{4(1-s)}{(s+0.1)(s+4)},
\end{equation}
and $C(s)=1$. Figure \ref{fig_ex2} shows the Bode magnitude and
phase plot of $P(s)$, $P(s)/Q_{1,2}(s)$, and $P(s)/Q_{1,4}(s)$. It
can be easily verified that in this example $P(s)$,
$P(s)/Q_{1,2}(s)$, and $P(s)/Q_{1,4}(s)$ lead to the phase margins
$2.7^\circ$, $32.7^\circ$, and $50.3^\circ$, and the gain margins
$0.2149$ dB, $4.401$ dB, and $8.958$ dB, respectively. As it is
observed the proposed cancellation method has effectively
increased the phase and the gain margin.

\textbf{Example 3.} In this example we show that the proposed
fractional-order unstable pole-zero cancellation strategy can
enhance the performance of the feedback system by simultaneous
decrement of overshoots, undershoots and the control effort. In
order to make a fair comparison, we study the behavior of feedback
system with and without using the proposed cancellation method,
where in both cases $C(s)$ in Fig. \ref{fig_prop} is calculated
optimally by using the LQG servo-controller design algorithm.

\begin{figure}
\begin{center}
%[height=4cm]
\includegraphics[width=8.4cm]{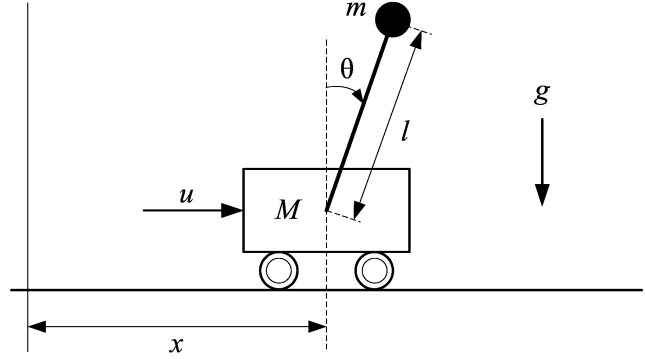}    % The printed column
\caption{Inverted pendulum, corresponding to Example 3}
\label{fig_pen}                                 % Size the figures
\end{center}                                 % accordingly.
\end{figure}

Linearizing the nonlinear equations governing the inverted
pendulum shown in Fig. \ref{fig_pen} around the unstable
equilibrium point at $(x,\theta)=(0,0)$ leads to the following
unstable non-minimum phase transfer function \cite{hoagg}:
\begin{equation}\label{p2}
P(s)=\frac{X(s)}{U(s)}=\frac{1}{M}\frac{(s-z)(s+z)}{s^2(s-p)(s+p)},
\end{equation}
where
\begin{equation}
p=\sqrt{\frac{g}{l}+\frac{mg}{Ml}},\quad z=\sqrt{\frac{g}{l}}.
\end{equation}

In the following we design an LQG servo-controller for the
linearized process by using and without using the proposed
fractional-order unstable pole-zero canceller assuming that
$M=m=0.1$ kg, $l=1$ m, and $g=9.8$ $\mathrm{m/s^2}$.

First, consider the feedback system shown in Fig.
\ref{fig_system}. The Matlab function \textsf{lqg} can be used to
design a controller for this system such that the cost function:
\begin{equation}
J=E\left\{\lim_{T\rightarrow\infty} \frac{1}{T}\int_{0}^T
\left([\mathbf{x}^\mathrm{T},u]\mathbf{R}\left[\begin{array}{c}
  \mathbf{x} \\
  u \\
\end{array}\right]+q e^2\right)\mathrm{d}t\right\},
\end{equation}
is minimized, where $\mathbf{x}$ is the state vector of process,
$u$ is the control applied to process, and $e$ is the tracking
error defined as $r-(y+n)$. Without any loss of generality, it is
assumed that in the problem under consideration $\mathbf{R}$ is a
diagonal matrix of suitable dimension where all non-zero entries
are equal to unity, and $q=1$. Moreover, it is assumed that the
process noise $w$ (which is not shown in Fig. \ref{fig_system})
and measurement noise $n$ are Gaussian white noises with
covariance
\begin{equation}
E\left\{\left[%
\begin{array}{c}
  \mathbf{w} \\
  n \\
\end{array}%
\right][\mathbf{w}^\mathrm{T},n]\right\}=\mathbf{M},
\end{equation}
where, again without the loss of generality, it is assumed that
$\mathbf{M}=\mathrm{diag} [0.1,\ldots, 0.1]$. In dealing with the
feedback system shown in Fig. \ref{fig_prop} the LQG
servo-controller can be designed in the same manner assuming that
the transfer function of process is equal to $P(s)/Q_{\lambda,v}
(s)$.

\begin{figure}
\begin{center}
%[height=4cm]
\includegraphics[width=8.4cm]{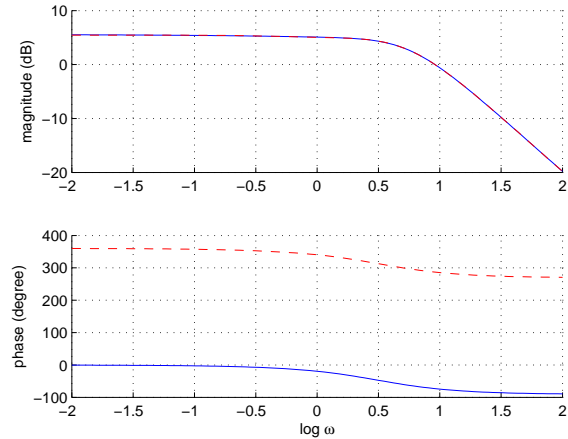}    % The printed column
\caption{Bode phase and magnitude plot of (\ref{ex2}) (solid
curve) and its approximation (dashed curve), corresponding to
Example 3}
\label{fig1}                                 % Size the figures
\end{center}                                 % accordingly.
\end{figure}

The simplest possible approach to control a process with the
transfer function given in (\ref{p2}) using the feedback system
shown in Fig. \ref{fig_prop} is to consider $Q_{\lambda,v}(s)$
equal to $s^{1/2}+z^{1/2}$. It can be shown that although this
choice can considerably decrease the overshoots and undershoots in
the response of feedback system to step command, it will apply
more control effort compared to the case such a cancellation is
not applied. To remove this difficulty, the fractional-order
pole-zero canceller in Fig. \ref{fig_prop} (i.e., the box with
transfer function $1/Q_{\lambda,v}(s)$) can be substituted with:
\begin{equation}\label{canc3}
\frac{Q_{p,2}(s)}{Q_{z,2}(s)}=\frac{s^{1/2}+p^{1/2}}{s^{1/2}+z^{1/2}}.
\end{equation}
Note that using the above canceller is equivalent to the half
cancellation of both unstable pole and non-minimum phase zero.
Now, according to the previous discussions, the LQG
servo-controller should be designed for a system with the
open-loop transfer function:
\begin{equation}\label{ex1}
\frac{(s^{1/2}-z^{1/2})(s+z)}{Ms^2\left(s^{1/2}-p^{1/2}\right)(s+p)}.
\end{equation}
For this purpose first we should approximate the transfer function
given in (\ref{ex1}) with an integer-order transfer function by
using an order-reduction algorithm. The main reason for working
with such an approximate integer-order transfer function is that
the standard LQG design algorithm can be applied only to linear
processes modelled by integer-order transfer functions. In order
to approximate (\ref{ex1}) with an integer-order transfer function
first we apply the Matlab function \textsf{invfreqs} to the
fractional-order transfer function
\begin{equation}\label{ex2}
\frac{s^{1/2}-z^{1/2}}{M\left(s^{1/2}-p^{1/2}\right)(s+p)},
\end{equation}
and then we add a zero (located at $s=-z$) and two poles (located
at $s=0$) to the resulted transfer function. Clearly, it is also
possible to apply the \textsf{invfreqs} command directly to
(\ref{ex1}), but that would lead to less accurate results compared
to the proposed approach. The solid and dashed curve in Fig.
\ref{fig1} show the Bode plots of (\ref{ex2}) and its
approximation, respectively (in this figure the approximation is
performed in the frequency range $[10^{-2},10^{2}]$ rad/s, where
the degree of the numerator and denominator of the approximating
transfer function is considered equal to 6). As it can be observed
in Fig. \ref{fig1} the two plots are in a fair agreement.

% and consequently, the LQG
%servo-controller can be designed for the corresponding
%integer-order approximate transfer function.

\begin{figure}
\begin{center}
%[height=4cm]
\includegraphics[width=8.4cm]{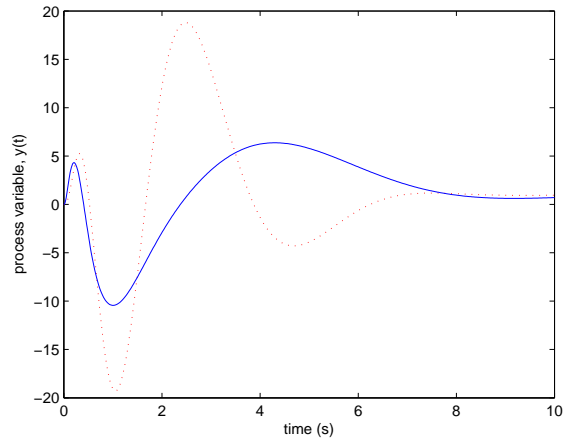}    % The printed column
\caption{Position of the LQG-controlled cart with (solid curve)
and without (dotted curve) using the fractional-order unstable
pole-zero canceller given in (\ref{canc3}), corresponding to
Example 3}
\label{fig3}                                 % Size the figures
\end{center}                                 % accordingly.
\end{figure}

The LQG servo-controller is designed for $P(s)$ and $Q_{p,2}(s)
P(s)/Q_{z,2}(s)$ (as described above) and the responses of
corresponding feedback systems to the unit step command are shown
in Fig. \ref{fig3}. In this figure the solid curve shows the cart
position when the fractional-order unstable pole-zero canceller
given in (\ref{canc3}) is applied, and the dotted curve shows the
cart position without using such a canceller. This figure clearly
shows that the proposed method can effectively improve the
response of feedback system by decreasing the overshoots and
undershoots.

\begin{figure}
\begin{center}
%[height=4cm]
\includegraphics[width=8.4cm]{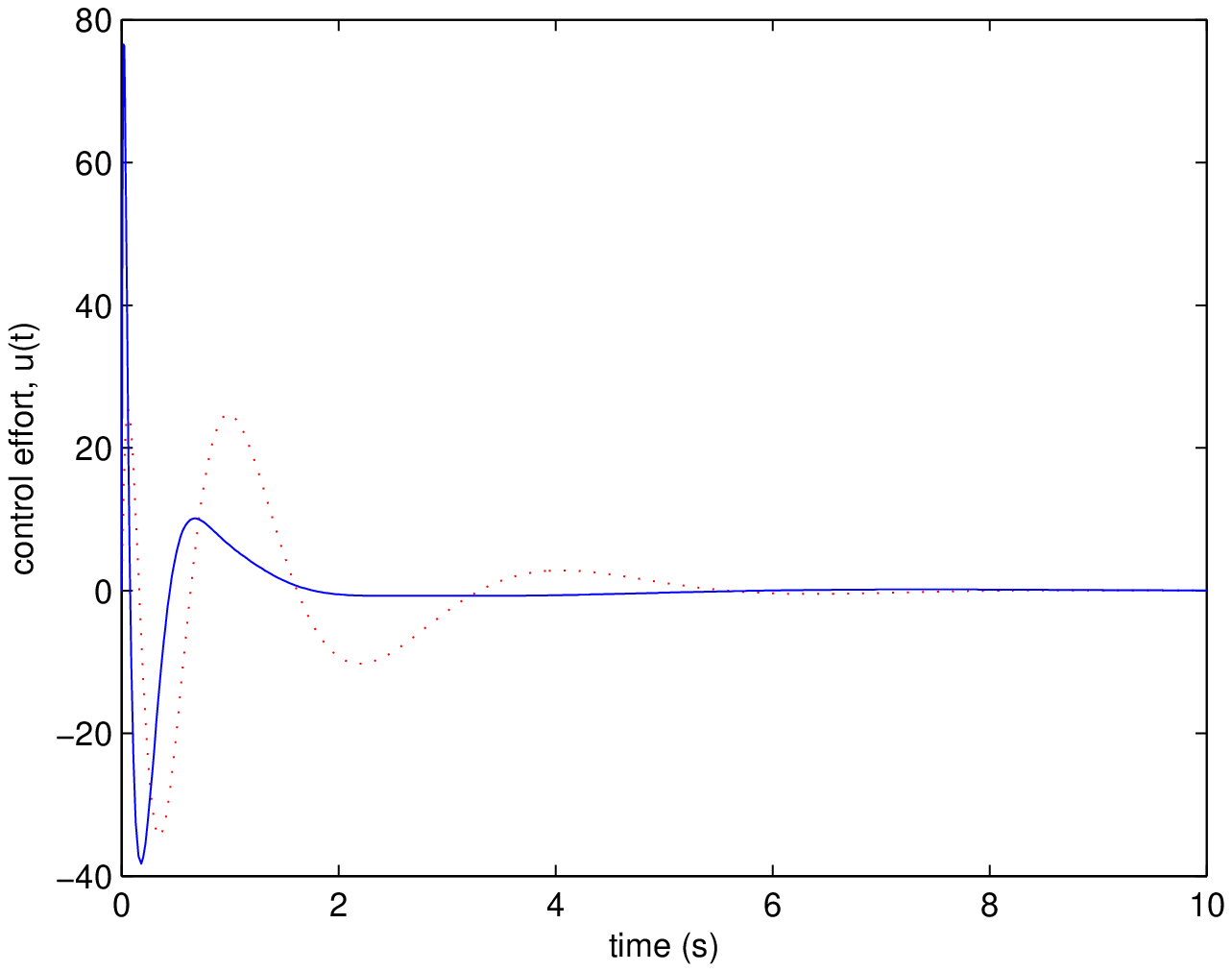}    % The printed column
\caption{Control effort with (solid curve) and without (dotted
curve) using the fractional-order unstable pole-zero canceller
given in (\ref{canc3}), corresponding to Example 3}
\label{fig4}                                 % Size the figures
\end{center}                                 % accordingly.
\end{figure}

Note that in order to simulate the feedback system at the presence
of fractional-order unstable pole-zero canceller given in
(\ref{canc3}) we have approximated the canceller with an
integer-order transfer function and then we have used it in the
connection shown in Fig. \ref{fig_prop} (in this example the
approximation is performed in the frequency range $[10^{-3},10^4]$
rad/s and the order of the numerator and denominator of the
approximating transfer function is considered equal to 4). Figure
\ref{fig4} shows the control effort with and without using the
fractional-order pole-zero canceller given in (\ref{canc3}).
Interesting observation is that applying this fractional-order
pole-zero canceller also decreases the control effort. More
precisely, energy of the control signal with and without using the
fractional-order unstable pole-zero canceller is equal to $483.5$
and $693.5$.

%Figure \ref{fig2} shows the Bode plots of (\ref{canc3}) and its
%integer-order approximation in the frequency range
%$[10^{-3},10^4]$ rad/s when the order of the numerator and
%denominator of the approximating transfer function is considered
%equal to 4.

%Numerical simulations performed by the author show that better
%results, i.e. better responses with smaller control efforts, can
%be obtained by better cancellation of the unstable pole and
%non-minimum phase zero of the process. It can be achieved by using
%e.g.
%\begin{equation}\label{canc4}
%\frac{(s^{1/2}+p^{1/2})(s^{1/4}+p^{1/4})}{(s^{1/2}+z^{1/2})(s^{1/2}+z^{1/4})},
%\end{equation}
%or
%\begin{equation}\label{canc5}
%\frac{(s^{1/2}+p^{1/2})(s^{1/4}+p^{1/4})(s^{1/8}+p^{1/8})}{(s^{1/2}+z^{1/2})(s^{1/2}+z^{1/4})(s^{1/8}+z^{1/8})},
%\end{equation}
%instead of (\ref{canc3}).

%\begin{figure}
%\begin{center}
%%[height=4cm]
%\includegraphics[width=8.4cm]{fig1.eps}    % The printed column
%\caption{Bode phase and magnitude plot of (\ref{canc3}) (solid
%curve) and its integer-order approximation (dashed curve)}
%\label{fig2}                                 % Size the figures
%\end{center}                                 % accordingly.
%\end{figure}

\section{Conclusion and discussion}\label{sec_conc}
In this paper we showed that although the unstable pole-zero
cancellation between controller and process is impractical and
leads to internal stability, the fractional-order pole-zero
cancellation is possible and can be effective. A method for
designing fractional-order unstable pole-zero cancellers is
proposed and it is specially shown that the proposed method can
increase the phase and the gain margin, and consequently, improve
the efficiency of the classical controller design algorithm under
consideration. Effect of the proposed method on time-domain
responses of the non-minimum phase open-loop system is studied and
it is shown that it can decrease the overshoots and undershoots in
its step response.

There are still many other questions to be answered. For example,
the fractional-order unstable pole-zero cancellation strategy can
be performed in many other ways. However, it is not exactly known
which one is more effective. Moreover, the controller can also be
designed using many different methods. It is expected that the
proposed fractional-order pole-zero cancellation method be more
effective when a certain controller design algorithm is used, but
it is not known which one is that. In fact, the proposed approach
may also be effective even when the process is minimum phase or
stable. These the questions that can be considered as the subject
of future studies.


\begin{thebibliography}{}

\end{thebibliography}


\begin{thebibliography}{99}     % Otherwise use the
                                 % thebibliography environment.
                                 % Insert the full references here.
                                 % See a recent issue of Automatica
                                 % for the style.
\bibitem{astrom}
K.J. \AA str\"{o}m, R.D. Bell, Drum-boiler dynamics,
\emph{Automatica}, 36 (2000) 363--378.

\bibitem{hauser}
J. Hauser, S. Sastry, G. Meyer, Nonlinear control design for
slightly non-minimum phase systems: Application to V/STOL
aircraft, \emph{Automatica}, 28(4) (1992) 665--679.

\bibitem{kwon}
D.-S. Kwon, W.J. Book, A time-domain inverse dynamic tracking
control of a single-link flexible manipulator, \emph{Journal of
Dynamic Systems, Measurement, and Control}, 116 (1994) 193--200.

\bibitem{kravaris}
C. Kravaris, P. Daoutidis, Nonlinear state feedback control of
second-order nonminimum-phase nonlinear systems, \emph{Computers
and Chemical Engineering}, 14(4/5) (1990) 439--449.

\bibitem{seron1}
M.M. Seron, J.H. Braslavsky, G.C. Goodwin, \emph{Fundamental
Limitations in Filtering and Control}, New York: Springer-Verlag,
1997.

\bibitem{seron2}
M.M. Seron, J.H. Braslavsky, P.V. Kokotovic, D.Q. Mayne, Feedback
limitations in nonlinear systems: From Bode integrals to cheap
control, \emph{IEEE Trans. Automat. Contr.} 44 (1999) 829--833.

\bibitem{qiu}
L. Qiu, E.J. Davison, Performance limitations of non-minimum phase
systems in the servomechanism problem, \emph{Automatica}, 29
(1993) 337--349.

\bibitem{middleton}
R.H. Middleton, Tradeoffs in linear control system design,
\emph{Automatica}, 27(2) (1991) 281--292.

\bibitem{hoagg}
J.B. Hoagg, D.S. Bernstein, Nonminimum-phase zeros, \emph{IEEE
Control Systems Magazine}, (2007) 45--57.

\bibitem{doyle}
J.C. Doyle, B.A. Francis, A.R. Tannenbaum, \emph{Feedback Control
Theory}, New York: Macmillan, 1992.

\bibitem{kailath}
T. Kailath, \emph{Linear Systems}, Englewood Cliffs, NJ:
Prentice-Hall, 1980.

\bibitem{Aguiar}
A. Pedro Aguiar, J.P. Hespanha, P.V. Kokotovic, Path-following for
nonminimum phase systems removes performance limitations,
\emph{IEEE Trans. Automat. Contr.} 50(2) (2005) 234--239.

\bibitem{fliness}
M. Fliess, H. Sira-Ramirez, R. Marquez, Regulation of nonminimum
phase outputs: A flatness based approach, in: D. Normand-Cyrot
(Ed.), \emph{Perspectives in Control-Theory and Applications: A
Tribute to Ioan Dor\'{e} Landau},  London, U.K.: Springer-Verlag,
1998, pp. 143–-164.

\bibitem{hiddabi}
S. Al-Hiddabi, N. McClamroch, Tracking and maneuver regulation
control for nonlinear nonminimum phase systems: Application to
flight control, \emph{IEEE Trans. Control Syst. Technol.} 10(6)
(2002) 780-–792.


\bibitem{silverman}
R.A. Silverman, \emph{Complex Analysis with Applications}, Dover
Publications, Inc., New York, 1984.


\bibitem{lau}
K. Lau, R.H. Middleton, J.H. Braslavsky, Undershoot and settling
time tradeoffs for nonminimum phase systems, \emph{IEEE Trans.
Automat. Contr.}, 48(8) (2003) 1389--1393.

\bibitem{matignon1}
D. Matignon, Stability properties for generalized fractional
differential systems, ESAIM: Proc., 1998, Vol. 5, pp. 145--158.


%\bibitem{matignon2}
%Matignon, D. (??????). Stability results for fractional fractional
%differential equations with application to control processing, ???











%
%
%\bibitem{farshad1}
%Merrikh-Bayat, F., \& Karimi-Ghartemani, M. (2008). On the
%essential instabilities caused by fractional-order transfer
%functions. \emph{Mathematical Problems in Engineering}, Volume
%2008, Article ID 419046, 13 pages, doi:10.1155/2008/419046.
%
%\bibitem{podlubny99}
%Podlubny, I. (1999). {\it Fractional differential equations}.
%Academic Press.



%
%\bibitem[Heritage, 1992]{Heritage:92}
%     (1992) {\it The American Heritage.
%     Dictionary of the American Language.}
%     Houghton Mifflin Company.
%  \bibitem[Able, 1956]{Abl:56}
%     B.~C.~Able (1956). Nucleic acid content of macroscope.
%     {\it Nature 2}, 7--9.
%  \bibitem[Able {\em et al.}, 1954]{AbTaRu:54}
%     B.~C. Able, R.~A. Tagg, and M.~Rush (1954).
%     Enzyme-catalyzed cellular transanimations.
%     In A.~F.~Round, editor,
%     {\it Advances in Enzymology Vol. 2} (125--247).
%     New York, Academic Press.
%  \bibitem[R.~Keohane, 1958]{Keo:58}
%     R.~Keohane (1958).
%     {\it Power and Interdependence:
%     World Politics in Transition.}
%     Boston, Little, Brown \& Co.
%  \bibitem[Powers, 1985]{Pow:85}
%     T.~Powers (1985).
%     Is there a way out?
%     {\it Harpers, June 1985}, 35--47.
%
\end{thebibliography}
\end{document}